\documentclass[leqno,11pt]{article}
\usepackage{amssymb}
\usepackage{amsfonts}
\usepackage{amsmath}
\usepackage{amsxtra}
\usepackage[all]{xy}

\newcommand\uu{\mathtt{u}}

\def\LL{\mathtt{L}}
\def\MM{\mathtt{M}}
\def\NN{\mathtt{N}}
\def\kk{\mathtt{k}}
\def\rr{{\mathfrak r}}

\title{New structures on the tangent bundles\\ and tangent sphere bundles}
\author{Marian Ioan Munteanu
\thanks{ Beneficiary of a CNR-NATO Advanced Research Fellowship pos. 216.2167 Prot. n. 0015506. }}
\date{ }

\newtheorem{contor}{1}[section]

\newtheorem{theorem}[contor]{Theorem}
\newtheorem{lemma}[contor]{Lemma}
\newtheorem{proposition}[contor]{Proposition}
\newtheorem{remark}[contor]{Remark}
\newtheorem{example}[contor]{Example}
\newtheorem{corollary}[contor]{Corollary}

\def\proof{{\sc Proof.\ }}
\newcommand{\gata}{\hfill\hskip -1cm \rule{.5em}{.5em}}

\begin{document}
\maketitle

\begin{abstract}

\noindent
In this paper we study a Riemanian metric on the tangent bundle $T(M)$ of a Riemannian
manifold $M$ which generalizes Sasaki metric and Cheeger Gromoll metric and a compatible
almost complex structure which together with the metric confers to $T(M)$ a structure of
locally conformal almost K\"ahlerian manifold. This is the natural generalization
of the well known almost K\"ahlerian structure on $T(M)$.
We found conditions under which $T(M)$ is almost K\"ahlerian, locally conformal K\"ahlerian
or K\"ahlerian or when $T(M)$ has
constant sectional curvature or constant scalar curvature.  Then we will restrict to the unit
tangent bundle and we find an isometry with the tangent sphere bundle (not necessary unitary)
endowed with the restriction of the Sasaki metric from $T(M)$.
Moreover, we found that this map preserves also the natural almost contact structures obtained
from the almost Hermitian ambient structures on the unit tangent bundle and the tangent
sphere bundle, respectively.

\


\


\end{abstract}

\section{Introduction}
A Riemannian metric $g$ on a smooth manifold $M$ gives rise to
several Riemannian metrics on the tangent bundle $T(M)$ of $M$.
Maybe the best known example is the Sasaki metric $g_S$ introduced
in \cite{kn:Sas58}. Although the Sasaki metric is {\em naturally}
defined, it is very {\em rigid} in the following sense. For
example, O.Kowalski \cite{kn:Kow71} has shown that the tangent
bundle $T(M)$ with the Sasaki metric is never locally symmetric
unless the metric $g$ on the base manifold is flat. Then, E.Musso
\& F.Tricerri \cite{kn:MT88} have shown a more general result,
namely, the Sasaki metric has constant scalar curvature if and
only if $(M,g)$ is locally Euclidian. In the same paper, they have
given an explicit expression of a positive definite Riemannian
metric introduced by J.Cheeger and D.Gromoll in \cite{kn:CG72} and
called this metric {\em the Cheeger-Gromoll metric}. In
\cite{kn:Sek91} M.Sekizawa computed the Levi Civita connection,
the curvature tensor, the sectional curvatures and the scalar
curvature of this metric. These results are completed in 2002 by
S.Gudmundson and E.Kappos in \cite{kn:GK02}. They have also shown
that the scalar curvature of the Cheeger Gromoll metric is never
constant if the metric on the base manifold has constant sectional
curvature. Furthermore, M.T.K.Abbassi \& M.Sarih have proved that
$T(M)$ with the Cheeger Gromoll metric is never a space  of constant
sectional curvature (cf. \cite{kn:AS03}). A more general metric
is given by M.Anastasiei in \cite{kn:Ana99} which generalizes
both of the two metrics mentioned above in the following sense: it
preserves the orthogonality of the two distributions, on the horizontal
distribution it is the same as on the base manifold, and finally the Sasaki and
the Cheeger Gromoll metric can be obtained as particular cases of
this metric. A compatible almost complex structure is also introduced
and hence $T(M)$ becomes an locally conformal almost K\"aherian manifold.

V.Oproiu and his collaborators constructed a family of Riemannian metrics on
the tangent bundles of Riemannian manifolds which possess interesting
geometric properties (cf. \cite{kn:Opr99, kn:Opr199, kn:Opr01, kn:OP04}).
In particular, the scalar curvature of $T(M)$ can be constant also for
a non-flat base manifold with constant sectional curvature. Then
M.T.K.Abbassi \& M.Sarih proved in \cite{kn:AS05} that the considered
metrics by Oproiu form a particular subclass of the so-called $g$-{\em natural
metrics} on the tangent bundle (see also \cite{kn:Aba04,
kn:AK05, kn:AS05, kn:AS105, kn:AS205, kn:KS88}).

By thinking $T(M)$ as a vector bundle associated with $O(M)$ (the space of orthonormal frames
on $M$), namely $T(M)\equiv O(M)\times {\bf R}^n / O(n)$ (where the orthogonal group
$O(n)$ acts on the right on $O(M)$), Musso \& Tricerri construct natural
metrics on $T(M)$ (see \S4 in \cite{kn:MT88}). The idea is to consider $Q$
a symmetric, semi-positive definite tensor field of type $(2,0)$ and rank $2n$ on
$O(M)\times{\bf R}^n$. Assuming that $Q$ is {\em basic} for
$\psi:O(M)\times{\bf R}^n\longrightarrow T(M)$,
$({\mathbf{u}},\zeta)\mapsto\left(p,\zeta^i\uu_i\right)$,
where ${\mathbf{u}}=(p,\uu_1,\ldots,\uu_n)$ and $\zeta=(\zeta^1,\ldots,\zeta^n)$
(i.e. $Q$ is $O(n)$-invariant
and $Q(X,Y)=0$ for all $X$ tangent to a fiber of $\psi$) there is a unique Roiemannian
metric $g_Q$ on $T(M)$ such that $\psi^*g_Q=Q$.
In this paper we will show that the metric introduced in \cite{kn:Ana99}
can be construct by using the method of Musso and Tricerri and we study it.
Then we will give the conditions under which
$T(M)$ is locally conformal K\"ahlerian and respectively K\"ahlerian (Theorems 2.6 and 2.8).
These results extend the known result saying that $T(M)$ endowed with the Sasaki metric and
the canonical almost complex structure is K\"ahlerian if and only if the base manifold is
locally Euclidean.

Next we want to have constant sectional curvature and constant scalar curvature,
respectively on $T(M)$. With this end in view, we compute the Levi Civita connection, the
curvature tensor, the sectional curvature and the scalar curvature of this metric.
We found relations between the sectional curvature (resp. scalar curvature) on $T(M)$ and
the corresponding curvature on the base $M$. We give an example of metric on $T(M)$ of
Cheeger Gromoll type which is flat. (Recall the fact that Cheeger Gromoll metric can not
have constant sectional curvature.) See Proposition \ref{prop:2.17}.
We also obtain a locally conformal K\"ahler structure
(cf. Example \ref{ex:2.8}$/_2$) and a K\"ahler structure (cf. Remark \ref{rem:2.10}$/_3$)
on $T(M)$.
We give some examples of metrics on $T(M)$ (when $M$ is
a space form) having constant scalar curvature.
See Examples \ref{ex:2.22} and \ref{ex:2.23}.

In section 3 we restrict the structure on the unit tangent bundle, obtaining an almost
contact metric. We will show that the unit tangent bundle is isometric with a tangent
sphere bundle $T_r(M)$ (we find the radius $r$) endowed with the restriction of Sasaki
metric from $T(M)$ (see also \cite{kn:BV00}, Remark 4, p.88). Moreover, this map preserves the almost contact structures.
M.Sekizawa \& O.Kowalski have studied the geometry of the tangent sphere bundles
with arbitrary radii endowed with the induced Sasaki metric (see \cite{kn:KS00}).
They have also noticed that the unit tangent bundle equipped with the induced Cheeger
Gromoll metric is isometric to the tangent sphere bundle $T_\frac1{\sqrt{2}}(M)$,
of radius $\frac 1{\sqrt{2}}$ endowed with the metric induced by the Sasaki metric.
Some other generalizations concerning this fact are given in \cite{kn:AK05}.
In the end of the section we obtained some properties for $T_1(M)$ as contact manifold. Among
the results we state the following: {\em The contact metric structure on $T_1(M)$ is
$K$-contact if and only if the base manifold has positive constant sectional curvature.
In this case $T_1(M)$ becomes a Sasakian manifold.}

\section{The tangent bundle $T(M)$}
Let $(M,g)$ be a Riemannian manifold and let $\nabla$ be its Levi
Civita connection. Let $\tau:T(M)\longrightarrow M$ be the tangent
bundle. If $\uu\in T(M)$ it is well known the following decomposition of the
tangent space $T_\uu T(M)$ (in $\uu$ at $T(M)$)
$$
    T_\uu T(M)=V_\uu T(M)\oplus H_\uu T(M)
$$
where $V_\uu T(M)=\ker \tau_{*,\uu}$ is the vertical space and
$H_\uu T(M)$ is the horizontal space in $\uu$ obtained by using $\nabla$.
(A curve $\widetilde\gamma:I\longrightarrow T(M)$ , $t\mapsto(\gamma(t),V(t))$
is {\em horizontal} if the vector field $V(t)$ is parallel along
$\gamma=\widetilde\gamma\circ\tau$. A vector on $T(M)$ is {\em horizontal}
if it is tangent to an horizontal curve and {\em vertical} if it is tangent
to a fiber. Locally, if $(U,x^i)$, $i=1, \ldots, m$, $m=\dim M$,
is a local chart in $p\in M$, consider
$(\tau^{-1}(U), x^i, y^i)$ a local chart on $T(M)$. If $\Gamma_{ij}^k(x)$ are the
Christoffel symbols, then
$\delta_i=\frac\partial{\partial x^i}-\Gamma_{ij}^k(x)y^j\ \frac\partial{\partial y^k}$
in $\uu$, $i=1,\ldots,m$ span the space $H_\uu T(M)$, while
$\frac\partial{\partial y^i}$, $i=1,\ldots,m$ span the vertical space $V_\uu T(M)$.)
We have obtained the horizontal (vertical) distribution $HTM$ ($VTM$) and a direct sum decomposition
$$
   TTM=HTM\oplus VTM
$$
of the tangent bundle of $T(M)$.
If $X\in\chi(M)$, denote by $X^H$ (and $X^V$, respectively) the horizontal lift
(and the vertical lift, respectively) of $X$ to $T(M)$.

If $\uu\in T(M)$ then we consider the energy density in $\uu$ on $T(M)$,
namely $t=\frac12\ g_{\tau(\uu)}(\uu,\uu)$.

{\bf The Sasaki metric}
is defined uniquely by the following relations
\begin{equation}
\left\{
\begin{array}{l}
g_S(X^H,Y^H)=g_S(X^V,Y^V)=g(X,Y)\circ\tau\\
g_S(X^H,Y^V)=0,
\end{array}
\right.
\end{equation}
for each $X,Y\in\chi(M)$.

On $T(M)$ we an also define an almost complex structure $J_S$ by
\begin{equation}
J_SX^H=X^V,\quad J_SX^V=-X^H, \forall X\in\chi(M).
\end{equation}

It is known that $(T(M),J_S,g_S)$ is an almost K\"ahlerian manifold. Moreover, the integrability
of the almost complex structure $J_S$ implies that $(M,g)$ is locally flat (see e.g. \cite{kn:Bla02}).

{\bf The Cheeger-Gromoll metric}
on $T(M)$ is given by
\begin{equation}
\left\{
\begin{array}{l}
{g_{CG}}_{(p,\uu)}(X^H,Y^H)=g_p(X,Y),\
{g_{CG}}_{(p,\uu)}(X^H,Y^V)=0\\
\displaystyle
{g_{CG}}_{(p,\uu)}(X^V,Y^V)=\frac1{1+2t}\ \left(g_p(X,Y)+g_p(X,\uu)g_p(Y,\uu)\right)
\end{array}
\right.
\end{equation}
for any vectors $X$ and $Y$ tangent to $M$.

Since the almost complex structure $J_S$ is no longer compatible with the metric
$g_{CG}$, one defines on $T(M)$ another almost complex structure $J_{CG}$, compatible
with the Chegeer-Gromoll metric, by the formulas
\begin{equation}
\left\{
\begin{array}{l}
\displaystyle
J_{CG}X^H_{(p,\uu)}=\rr X^V-\frac1{1+\rr}\ g_p(X,\uu)\uu^V\\ \ \\
\displaystyle
J_{CG}X^V_{(p,\uu)}=-\frac1\rr\ X^H-\frac1{\rr(1+\rr)}\ g_p(X,\uu)\uu^H
\end{array}
\right.
\end{equation}
where $\rr=\sqrt{1+2t}$ and $X\in T_p(M)$. Remark that $J_{CG}\uu^H=\uu^V$
and $J_{CG}\uu^V=-\uu^H$. We get an almost Hermitian manifold $(T(M),J_{CG},g_{CG})$.
Moreover, if we denote by $\Omega_{CG}$ the Kaehler 2-form (namely
$\Omega_{CG}(U,V)=g_{CG}(U,J_{CG}V), \forall U,V\in\chi(T(M)))$ it is
quite easy to prove the following
\begin{proposition}
We have
\begin{equation}
d\Omega_{CG}=\omega\wedge\Omega_{CG},
\end{equation}
where $\omega\in\Lambda^1(T(M))$ is defined by $\omega_{(p,\uu)}(X^H)=0$ and
$\omega_{(p,\uu)}(X^V)=-\left(\frac1{\rr^2}+\frac1{1+\rr}\right)g_p(X,\uu)$,
$X\in T_p(M)$.
\end{proposition}
\proof
A simple computation gives

\qquad $\Omega_{CG}(X^H,Y^H)=\Omega(X^V,Y^V)=0$

\qquad $\Omega_{CG}(X^H,Y^V)=-\frac1\rr\left(g(X,Y)+\frac1{1+\rr}\ g(X,\uu)g(Y,\uu)\right)$.

(From now on we will omit the point $(p,\uu)$.)

The differential of $\Omega_{CG}$ is given by

\qquad $d\Omega_{CG}(X^H,Y^H,Z^H)=d\Omega_{CG}(X^H,Y^H,Z^V)=d\Omega_{CG}(X^V,Y^V,Z^V)=0$

\qquad $d\Omega_{CG}(X^H,Y^V,Z^V)=\frac1\rr\left(\frac1{\rr^2}+\frac1{1+\rr}\right)
        \left[g(X,Y)g(Z,\uu)-g(X,Z)g(Y,\uu)\right]$

for any $X,Y,Z\in\chi(M)$.

Hence the statement. \gata

\begin{remark}\rm
The almost Hermitian manifold $(T(M),J_{CG},g_{CG})$ is never almost Kaehlerian (i.e. $d\Omega_{CG}\neq0$).
\end{remark}

Finally, we obtain a necessary condition for the integrability of $J_{CG}$  namely, the base manifold
$(M,g)$ should be locally Euclidian.

\vspace{2mm}

{\bf A general metric}, let's call it $g_A$, is in fact a family of Riemannian
metrics (depending on two parameters) and the Sasaki metric and the Cheeger-Gromoll
metric are obtained by taking particular values for the two parameters. It is
defined (cf. \cite{kn:Ana99}) by the following formulas
\begin{equation}
\left\{
\begin{array}{l}
{g_A}_{(p,\uu)}(X^H,Y^H)=g_p(X,Y)\\
{g_A}_{(p,\uu)}(X^H,Y^V)=0\\
{g_A}_{(p,\uu)}(X^V,Y^V)=a(t)g_p(X,Y)+b(t)g_p(X,\uu)g_p(Y,\uu),
\end{array}
\right.
\end{equation}
for all $X,Y\in\chi(M)$, where $a,b:[0,+\infty)\longrightarrow[0,+\infty)$ and $a>0$.
For $a=1$ and $b=0$ one obtains the Sasaki metric and for $a=b=\frac1{1+2t}$ one gets
the Cheeger-Gromoll metric.

\begin{proposition}
The metric defined above can be construct by using the method described by Musso
and Tricerri in {\rm \cite{kn:MT88}}.
\end{proposition}
\proof
If we denote by $\theta=(\theta^1,\ldots,\theta^n)$ the canonical 1-form on $O(M)$ (namely, if
${\mathbf{p}}:O(M)\longrightarrow M$, $\theta$ is defined by
$d{\mathbf{p}}_{\mathbf{u}}(X)=\theta^i(X)\uu_i$, for
${\mathbf{u}}=(p,\uu_1,\ldots,\uu_n)$ and $X\in T_p(M)$)
we have $R^*_{\mathbf{u}}(\theta^i)=(a^{-1})^i_h\theta^h$ for each $a\in O(n)$.
The vertical distribution of $\psi$ is defined by
$$\theta^i=0,\ \ D\zeta^i:=d\zeta^i+\zeta^j\omega^i_j$$
where $\omega=(\omega^j_i)_{i,j}$ denotes the $so(n)$-valued
connection 1-form defined by the Levi Civita connection of $g$.
Since $R^*_a(\omega^i_j)=(a^{-1})^i_h\omega^h_ka^k_j$ we can also
write $R^*_a(D\zeta^i)=(a^{-1})^i_hD\zeta^h$, for all $a\in O(n)$.

Consider now the following bilinear form on $O(M)$
\begin{equation}
\begin{array}{l}
Q_A=\sum\limits_{i=1}^n(\theta^i)^2+a(\frac12||\zeta||^2)\sum\limits_{i=1}^n(D\zeta^i)^2+
     b(\frac12||\zeta||^2)\left(\sum\limits_{i=1}^n\zeta^iD\zeta^i\right)^2.
\end{array}
\end{equation}
It is symmetric, semi-positive definite and basic. Moreover, since the following
diagram
$$
 \xymatrix{
O(M)\times {\mathbf{R}}^n \ar[d]_{\rm{proj}_1} \ar[r]^-{\psi} & T(M)\ar[d]^{\tau}\\
O(M) \ar[r]^-{\mathbf{p}} & M
 }
$$
commutes, we have $\psi^*g_A=Q_A$. (See for details \S4 in
\cite{kn:MT88}.)

\gata

Again, we have to find an almost complex structure on $T(M)$, call it $J_A$, which is compatible
with the metric $g_A$. Inspired from the previous cases we look for the almost complex structure $J_A$
in the following way
\begin{equation}
\left\{
\begin{array}{l}
\displaystyle
J_A X^H_{(p,\uu)}=\alpha X^V + \beta g_p(X,\uu)\uu^V\\ \ \\
\displaystyle
J_A X^V_{(p,\uu)}= \gamma X^H + \rho g_p(X,\uu)\uu^H
\end{array}
\right.
\end{equation}
where $X\in\chi(M)$ and $\alpha$, $\beta$, $\gamma$ and $\rho$ are smooth functions on $T(M)$
which will be determined from $J_A^2=-I$ and from the compatibility conditions with the metric
$g_A$. Following the computations made in \cite{kn:Ana99} we get first $\alpha=\pm\frac1{\sqrt a}$
and $\gamma=\mp \sqrt{a}$. Without lost of the generality we can take
$$\alpha=\frac1{\sqrt a} {\quad\rm  and\quad} \gamma=- \sqrt{a}.$$
Then one obtains
$$
\begin{array}{l}
\displaystyle
\beta=-\frac1{2t}\left(\frac1{\sqrt{a}}+\epsilon\frac1{\sqrt{a+2bt}}\right) {\quad \rm and\quad}
\rho=\ \frac1{2t}\left(\sqrt{a}+\epsilon\sqrt{a+bt}\right)
\end{array}
$$
where $\epsilon=\pm1$.

We have the almost complex structure $J_A$
\begin{equation}
\left\{
\begin{array}{l}
\displaystyle
J_AX^H=\frac1{\sqrt{a}}X^V-\frac1{2t}\left(\frac1{\sqrt{a}}+\epsilon\frac1{\sqrt{a+2bt}}\right)g(X,\uu)\uu^V\\
\displaystyle
J_AX^V=-\sqrt{a}X^H+\frac1{2t}\left(\sqrt{a}+\epsilon\sqrt{a+2bt}\right)g(X,\uu)\uu^H
\end{array}
\right.
\end{equation}
and the almost Hermitian manifold $(T(M),g_A,J_A)$.
\begin{remark}\rm
In this general case $J_A$ is defined on $T(M)\setminus 0$ (the bundle of non zero tangent vectors),
but if we consider $\epsilon=-1$ the previous relations define $J_A$ on all $T(M)$.
\end{remark}
\begin{remark}\rm
If we take $\epsilon=-1$, $a=1$ and $b=0$ we get the manifold $(T(M),g_S,J_S)$
and for $\epsilon=-1$, $a=b=\frac1{1+2t}$ we obtain the manifold $(T(M),g_{CG},J_{CG})$.
\end{remark}

If we denote by $\Omega_A$ the K\"ahler 2-form (i.e. $\Omega_A(U,V)=g_A(U,J_AV), \forall U,V\in\chi(T(M))$)
one obtains
\begin{proposition}
{\rm (see \cite{kn:Ana99})}
The almost Hermitian manifold $(T(M),g_A,J_A)$ is locally conformal almost K\"ahlerian, that is
\begin{equation}
d\Omega_A=\omega\wedge\Omega_A
\end{equation}
where $\omega$ is a closed and globally defined $1-$form on $T(M)$ given by
$$\omega(X^H)=0 {\quad and \quad}
\omega(X^V)=\frac1{\sqrt{a}}\left(\frac{a'}{\sqrt{a}}+\frac1{2t}\ (\sqrt{a}+\epsilon\sqrt{a+2bt} )\right)g(X,\uu).$$

\end{proposition}

As consequence one can state the following

\begin{theorem}
\label{th:aK}
The almost Hermitian manifold $(T(M),g_A,J_A)$ is almost K\"ahlerian if and only if
$$
b(t)=\frac{2a'(t)\left(ta'(t)+a(t)\right)}{a(t)}
$$
and for $\epsilon=-1$, $a(t)$ is an increasing function, while for $\epsilon=+1$, $ta(t)$
is a decreasing function.
\end{theorem}
\proof
The condition $\omega=0$ is equivalent to
$$2ta'(t)+a(t)=-\epsilon\sqrt{a(t)}\cdot\sqrt{a(t)+2tb(t)}\ .$$
From here, we get $b(t)$. Moreover it follows $(a(t)\sqrt{t}$ is a
monotone function, namely it is increasing if $\epsilon = -1$ and
decreasing for $\epsilon = +1$. Since $b(t)$ is positive we conclude

\quad $\bullet$ if $\epsilon=-1$: $2a't+a>0\longleftrightarrow 2(a't+a)>a\longrightarrow a't+a>0$

\qquad $\longrightarrow a'>0\longrightarrow$ $a$ increases (this implies $a\sqrt{t}$, $at$ are also increasing functions);

\quad $\bullet$ if $\epsilon=+1$: $2a't+a<0\longleftrightarrow a't+a<-a't\longrightarrow a't+a<0$

\qquad $\longrightarrow$ $at$ decreases (this implies $a\sqrt{t}$, $a$ are also decreasing functions).

\gata

\

{\bf The integrability of $J_A$.}

In order to have an integrable structure $J_A$ on $T(M)$ we have to compute the Nijenhuis
tensor $N_{J_A}$ of $J_A$ and to ask that it vanishes identically.

For the integrability tensor $N_{J_A}$ we have the following
relations
\begin{equation}
\label{eq:NIJ}
\left\{
\begin{array}{l}
N_{J_A}(X^H,Y^H)=\left(-\frac{a'}{2a^2}+\frac{a+ta'}{a\sqrt{a}}\ A(t) \right)\big(g(X,\uu)Y-g(Y,\uu)X\big)^V+(R_{XY}\uu)^V\\
N_{J_A}(X^V,Y^V)=\big(-aR_{XY}\uu+\sqrt{a}\ B(t)g(Y,\uu)R_{X\uu}\uu-\sqrt{a}\ B(t)R_{Y\uu}\uu\big)^V-\\
\qquad\qquad -\frac1{\sqrt{a}}\left(\frac {a'}{2\sqrt{a}}+B(t)\right)\big(g(Y,\uu)X-g(X,\uu)Y\big)^V
\end{array}
\right.
\end{equation}
where $A(t)=\frac 1{2t}\left(\frac1{\sqrt{a}}+\epsilon \frac1{\sqrt{a+2bt}}\right)$ and
$B(t)=\frac 1{2t}\left(\sqrt{a}+\epsilon\sqrt{a+2bt}\right)$. (The expression for $N_{J_A}(X^H,Y^V)$
is very complicated.)

Thus if $J_A$ is integrable then
$$
R_{XY}\uu=\left(-\frac{a'}{2a^2}+\frac{a+ta'}{a\sqrt{a}}\ A(t) \right)\big(g(Y,\uu)X-g(X,\uu)Y\big)
$$
for every $X,Y\in\chi(M)$ and for every point $\uu\in T(M)$. It follows that $M$ is a space form $M(c)$
($c$ is the constant sectional curvature of $M$).
Consequently,
\begin{equation}
\label{eq:c}
-\frac{a'}{2a^2}+\frac{a+ta'}{a\sqrt{a}}\ A(t)=c.
\end{equation}
So

\qquad\checkmark\qquad given $a(t)$ and $c$ we can easily find $b(t)$;

\qquad\checkmark\qquad given $b(t)$ and $c$ we have to solve an ODE in order to find $a(t)$;

\qquad\checkmark\qquad given $a(t)$ and $b(t)$ we have to check if $c$ in (\ref{eq:c}) is constant.

\begin{example}\rm\
\label{ex:2.8}

{\bf 1.} In Sasaki case $(a(t)=1, b(t)=0, \epsilon=-1)$ it follows $c=0$ i.e. $M$ is flat.

{\bf 2.} Looking for a locally conformal K\"ahler structure on $T(M)$ with the metric having $a(t)=b(t)$ we obtain
$$
a(t)=b(t)=\frac{e^{2\sqrt{1+2t}}}{2\left(ce^{2\sqrt{1+2t}}t+(1+t+\sqrt{1+2t})k\right)}
$$
with $k$ a positive real constant and $c$ must be nonnegative.
\end{example}
Replacing the expression of the curvature $R$ in (\ref{eq:NIJ}$){}_2$ we obtain again (\ref{eq:c}).

\vspace{2mm}

{\bf Question:} {\em Can $(T(M),g_A,J_A)$ be a Kaehler manifold?}

If this happens then the base manifold is a space form $M(c)$
and the functions $a$ and $b$ satisfy
\begin{equation}
\label{eq:K}
b=\frac{2a'(ta'+a)}{a} \quad{\rm and }\quad
\end{equation}
\begin{equation}
\label{eq:ODE}
a'=2ca(2ta'+a).
\end{equation}
If $c=0$ ($M$ is flat) then $a$ is a positive constant and $b$ vanishes.
\par
If $c\neq0$ the ODE (\ref{eq:ODE}) has general solutions
\begin{equation}
\label{eq:solODE}
a_{1,2}(t)=\frac{1\pm\sqrt{1+\kappa t}}{4ct}
\end{equation}
with $\kappa$ a real constant. Taking into account that $a$ and $b$ are
positive functions, using (\ref{eq:K}) one gets:

\textsc{Case 1.}

\begin{equation}
\label{eq:ab_case1}
a=\frac{1+\sqrt{1+\kappa t}}{4ct} \ {\rm and}\
b=-\frac{\kappa (1+\sqrt{1+\kappa t})}{8ct(1+\kappa t)}.
\end{equation}
\qquad
Here $c>0$, $t>0$, $\kappa<0$, $t<-\frac1\kappa$ and $\epsilon=+1$.

\textsc{Case 2.}

\begin{equation}
\label{eq:ab_case2}
a=-\frac{\kappa}{4c(1+\sqrt{1+\kappa t})} \ {\rm and}\
b=\frac{\kappa^2}{8c(1+\kappa t)(1+\sqrt{1+\kappa t})}.
\end{equation}

\qquad
Here $\kappa c<0$, $c<0$ (then $\kappa>0$), $t<-\frac1\kappa$ and $\epsilon=-1$.

Consider
$B_\kappa=\left\{\texttt{v}\in T(M)\ :
    \ g_{\tau(\texttt{v})}(\texttt{v},\texttt{v})<-\frac2\kappa\right\}$
    and
$\dot B_\kappa=B_\kappa\setminus M$.

\begin{theorem}
The manifolds $B_\kappa$ in {\textsc{Case 1}} and $\dot B_\kappa$ in {\textsc{Case 2}}
are Kaehler manifolds.
\end{theorem}

\begin{remark}\rm
\label{rem:2.10}
In order to have a positive definite metric $g_A$, the necessary and sufficient
conditions are $a>0$ and $a+2bt>0$ ($b>0$ is too strong).
Hence, the previous theorem can be reformulated as:

1.
$(T(M)\setminus M, g_A, J_A)$ where $a$ and $b$ are given by (\ref{eq:ab_case1}),
$c>0$ and $\epsilon=+1$ is a Kaehler manifold.

2.
$(B_\kappa, g_A, J_A)$ where $a$ and $b$ are given by (\ref{eq:ab_case2}),
$c>0$, $k<0$ and $\epsilon=-1$ is a Kaehler manifold.

3.
$(T(M),g_A,J_A)$ where $a$ and $b$ are given by (\ref{eq:ab_case2}),
$c<0$, $k>0$ and $\epsilon=-1$ is a Kaehler manifold.
\end{remark}

Now we give

\begin{proposition}
Let $(M,g)$ be a Riemannian manifold and let $T(M)$ be its tangent bundle equipped with the metric $g_A$.
Then, the corresponding Levi Civita connection $\tilde\nabla^A$ satisfies the following relations:
\begin{equation}
\left\{
\begin{array}{l}
\tilde\nabla^A_{X^H}Y^H=(\nabla_XY)^H-\frac12\ (R_{XY}\uu)^V\\
\ \\
\tilde\nabla^A_{X^H}Y^V=(\nabla_XY)^V+\frac a2\ (R_{\uu Y}X)^H\\
\ \\
\tilde\nabla^A_{X^V}Y^H=\frac a2\ (R_{\uu X}Y)^H\\
\ \\
\tilde\nabla^A_{X^V}Y^V=\LL \left(g(X,\uu)Y^V+g(Y,\uu)X^V\right)+\MM g(X,Y)\uu^V+
        \NN g(X,\uu)g(Y,\uu)\uu^V,
\end{array}
\right.
\end{equation}
where $\LL=\frac{a'(t)}{2a(t)}$, $\MM=\frac{2b(t)-a'(t)}{2(a(t)+2tb(t))}$
and $\NN=\frac{a(t)b'(t)-2a'(t)b(t)}{2a(t)(a(t)+2tb(t))}$.
\end{proposition}
\proof
The statement follows from Koszul formula making usual computations.
\gata

Having determined Levi Civita connection, we can compute now the Riemannian
curvature tensor $\tilde{R}^A$ on $T(M)$. We give

\begin{proposition}
The curvature tensor is given by
\begin{equation}
\label{eq:tildeRA}
\left\{
\begin{array}{l}
\tilde R^A_{X^HY^H}Z^H=(R_{XY}Z)^H+\frac a4\ \left[
   R_{\uu R_{XZ}\uu}Y-R_{\uu R_{YZ}\uu}X+2R_{\uu R_{XY}\uu}Z
    \right]^H+\\
    \qquad\qquad\qquad +\frac 12\ \left[
   (\nabla_ZR)_{XY}\uu
    \right]^V\\
\ \\
\tilde R^A_{X^HY^H}Z^V=\left[ R_{XY}Z+
  \frac a4\ (R_{YR_{\uu Z}X}\uu-R_{XR_{\uu Z}Y}\uu)
    \right]^V+\LL g(Z,\uu)(R_{XY}\uu)^V+\\
   \qquad\qquad\qquad +\MM g(R_{XY}\uu,Z)\uu^V+\frac a2\ \left[
(\nabla_XR)_{\uu Z}Y-(\nabla_YR)_{\uu Z}X
    \right]^H\\
\ \\
\tilde R^A_{X^HY^V}Z^H=\frac a2\left[(\nabla_XR)_{\uu Y}Z\right]^H+\\
\qquad\qquad\qquad
+\frac12\left[
    R_{XZ}Y-\frac a2\ R_{XR_{\uu Y}Z}\uu+\LL g(Y,\uu)R_{XZ}\uu+\MM g(R_{XZ}\uu,Y)\uu
    \right]^V\\
    \ \\
\tilde R^A_{X^HY^V}Z^V=-\frac a2\ (R_{YZ}X)^H-\frac{a^2}4\ (R_{\uu Y}R_{\uu Z}X)^H+\\
\qquad\qquad\qquad
 +\frac{a'}4\left[ g(Z,\uu)(R_{\uu Y}X)^H-g(Y,\uu)(R_{\uu Z}X)^H
 \right]\\
 \ \\
\tilde R^A_{X^VY^V}Z^H=a (R_{XY}Z)^H+\frac{a'}2\left[
  g(X,\uu)R_{\uu Y}Z-g(Y,\uu)R_{\uu X}Z
    \right]^H+\\
    \qquad\qquad\qquad+\frac {a^2}4\left[
 R_{\uu X}R_{\uu Y}Z-R_{\uu Y}R_{\uu X}Z
\right]^H\\
\ \\
\tilde R^A_{X^VY^V}Z^V=F_1(t) g(Z,\uu)\left[g(X,\uu)Y^V-g(Y,\uu)X^V\right]+\\
      \qquad\qquad\qquad +F_2(t)\left[g(X,Z)Y^V-g(Y,Z)X^V\right]+\\
      \qquad\qquad\qquad + F_3(t)\left[g(X,Z)g(Y,\uu)-g(Y,Z)g(X,\uu)\right]\uu^V,
\end{array}
\right.
\end{equation}
where $F_1=\LL'-\LL^2-\NN(1+2t\LL)$, $F_2=\LL-\MM(1+2t\LL)$ and $F_3=\NN-(\MM'+\MM^2+2t\MM\NN)$.
\end{proposition}

\begin{remark}\rm \

 (a) In the case of Sasaki metric we have:
$$
\LL=\MM=\NN=0\ ,\ F_1=F_2=F_3=0.
$$
 (b) In the case of Cheeger Gromoll metric we have (see also \cite{kn:GK02, kn:Sek91}):
$$
\begin{array}{l}
\LL=-\frac 1\rr\ ,\ \MM=\frac{\rr+1}{\rr^2}\ ,\ \NN=\frac 1{\rr^2}\ ,\
\LL'=\frac2{\rr^2}\ ,\ \MM'=-\frac{2(\rr+2)}{\rr^3}\ ,\ 1+2t\LL=\frac1\rr\\
F_1=\frac{\rr-1}{\rr^3}\ ,\ F_2=-\frac{\rr^2+\rr+1}{\rr^3}\ ,\ F_3=\frac{\rr+2}{\rr^3}
\end{array}
$$
where $\rr=1+2t$.
\end{remark}

In the following let $\tilde Q^A(U,V)$ denote the square of the area of the parallelogram with
sides $U$ and $V$ for $U,V\in\chi(T(M))$,
$$\tilde Q^A(U,V)=g_A(U,U)g_A(V,V)-g_A(U,V)^2.
$$
We have
\begin{lemma}
Let $X,Y\in T_pM$ be two orthonormal vectors. Then
\begin{equation}
\left\{
\begin{array}{l}
\tilde Q^A(X^H,Y^H)=1\\
\tilde Q^A(X^H,Y^V)=a(t)+b(t)g(Y,\uu)^2\\
\tilde Q^A(X^V,Y^V)=a(t)^2+a(t)b(t)\left(g(X,\uu)^2+g(Y,\uu)\right).
\end{array}
\right.
\end{equation}
\end{lemma}

We compute now the sectional curvature of the Riemannian manifold $(T(M), g_A)$, namely
$$\tilde K^A(U,V)=\frac{g_A(\tilde R^A_{UV}V,U)}{\tilde Q^A(U,V)}$$
for $U,V\in\chi(T(M))$.
\begin{proposition}
If $X,Y\in T_pM$ are two orthonormal vectors, then
\begin{equation}
\label{eq:tildeKA}
\left\{
\begin{array}{l}
\tilde K^A(X^H,Y^H)=K(X,Y)-\frac{3a(t)}4\ \left|R_{XY}\uu\right|^2\\
\tilde K^A(X^H,Y^V)=\frac{a(t)^2}{4(a(t)+b(t)g(Y,\uu)^2}\ \left|R_{\uu Y}X\right|^2\\
\tilde K^A(X^V,Y^V)=-\frac{F_1(t)a(t)g(Y,\uu)^2+F_2(t)\left(a(t)+b(t)g(X,\uu)^2\right)+F_3(t)(a(t)+2tb(t))g(X,\uu)^2}
    {a(t)^2+a(t)b(t)\left(g(X,\uu)^2+g(Y,\uu)^2\right)}\ .
\end{array}
\right.
\end{equation}
where $K(X,Y)$ is the sectional curvature of the plane spanned by $X$ and $Y$.
\rm
Here $|\cdot|$ denotes the norm of the vector with respect to the metric $g$ (in a point).
\end{proposition}
\proof
Calculations using (\ref{eq:tildeRA}) show that

\qquad
$g_A(\tilde R^A_{X^HY^H}Y^H,X^H)=g(R_{XY}Y,X)-\frac{3a}4\left|R_{XY}\uu\right|^2$

\qquad
$g_A(\tilde R^A_{X^HY^V}Y^V,X^H)=\frac{a^2}{4}\ \left|R_{\uu Y}X\right|^2$

\qquad
$g_A(\tilde R^A_{X^VY^V}Y^V,X^V)=-aF_1(t)g(Y,\uu)^2-F_2(t)(a+bg(X,\uu)^2)-
   F_3(t)g(X,\uu)^2(a+2tb)$.

Hence the conclusion.
\gata

Moreover, if $M$ has constant sectional curvature $c$, then $|R_{XY}\uu|^2=c^2\left(g(X,\uu)^2+g(Y,\uu)^2\right)$
and $|R_{\uu Y}X|^2=c^2g(\uu,X)^2$ for any orthonormal $X,Y\in T_p(M)$. Then we have

\qquad
$\begin{array}{l}
\tilde K^A(X^H,Y^H)=c-\frac{3a(t)c^2}4\ \left(g(X,\uu)^2+g(Y,\uu)^2\right)\\
\tilde K^A(X^H,Y^V)=\frac{a(t)^2c^2g(\uu,X)^2}{4(a(t)+b(t)g(Y,\uu)^2}\ .
\end{array}$

Following an idea from \cite{kn:Sek91} we are interested to study the sign of $K^A$. We have

{\bf (1)} If $c<0$ then $\tilde K^A(X^H,Y^H)<0$, if $c=0$ then $\tilde K^A(X^H,Y^H)=0$ and if $c>0$ then

\qquad
$\left\{
\begin{array}{l}
\tilde K^A(X^H,Y^H)>0,\ {\rm for\ } c\in(0,{\texttt{C}}_1)\\
\tilde K^A(X^H,Y^H)=0,\ {\rm for\ } c=\texttt{C}_1\\
\tilde K^A(X^H,Y^H)<0,\ {\rm for\ } c>\texttt{C}_1
\end{array}\right.
$, where $\texttt{C}_1=\frac 4{3a(t)\left(g(X,\uu)^2+g(Y,\uu)^2\right)}$.

Moreover the maximum value for $\tilde K^A(X^H,Y^H)$ is $\tilde K^A_{\rm max}=\frac {\texttt{C}_1}4$.

It will be better to have a constant $\texttt{C}>0$ (which does not depend on $X,Y$ and $t$) in the place of
$\texttt{C}_1$ so, we are looking for $\texttt{C}<\frac 4{3a(t)\left(g(X,\uu)^2+g(Y,\uu)^2\right)}$ for all
$X,Y$ and for any point $\uu$ of $T(M)$. We know that $g(X,\uu)^2+g(Y,\uu)^2\leq 2t$ for any $X,Y$ orthonormal,
so, to have this, it is sufficient for the function $a$ to verify
$$a(t)\leq\frac2{3\texttt{C}t}
$$
for any $t>0$. Remark that in the case of Cheeger Gromoll metric this fact occurs with $\texttt{C}=\frac43$
which is the best constant.

Another remark is that the maximum value for $\tilde K^A(X^H,Y^H)$ can be attained: for example
take $X=v$ with $|v|=1$, $Y=\frac1{\sqrt{2t-g(u,v)^2}}\ (u-g(u,v)v)$, $a=\frac2{3\texttt{C}t+o}$,
with $o$ a positive constant. If we take $\uu$ such that $t=\frac o3$ and $c=\frac{1+\texttt{C}}2$
one gets a maximum value $\frac{1+\texttt{C}}4$.

{\bf (2)} $\tilde K^A(X^H,Y^V)$ is non negative for all $c\neq0$ and vanishes for $c=0$.

{\bf (3)} Unfortunately for $\tilde K^A(X^V,Y^V)$ we cannot say anything yet.

\vspace{2mm}

Denote by $T_0(M)=T(M)\setminus{\mathbf{0}}$ the tangent bundle of non-zero vectors tangent to $M$. For a given
point $(p,\uu)\in T_0(M)$ consider an orthonormal basis $\{e_i\}_{i=\overline{1,m}}$ for the tangent space
$T_p(M)$ of $M$ such that $e_1=\frac u{|u|}$. Consider on $T_{(p,\uu)}T(M)$ the following vectors

\begin{equation}
\label{eq:E}
\left\{
\begin{array}{l}
E_i=e_i^H,\ i=1,\ldots,m\\
\ \\
E_{m+1}=\frac1{\sqrt{a+2tb}}\ e_1^V\\
\ \\
E_{m+k}=\frac 1{\sqrt{a}}\ e_k^V,\ k=2,\ldots,m.
\end{array}
\right.
\end{equation}
It is easy to check that $\{E_1,\ldots,E_{2m}\}$ is an orthonormal basis in $T_{(p,\uu)}T(M)$ (with respect to the metric $g_A$).
We will rewrite the expressions of the sectional curvature $\tilde K^A$ in terms of this basis. We have
\begin{equation}
\left\{
\begin{array}{ll}
\tilde K^A(E_i,E_j)=K(e_i,e_j)-\frac {3a(t)}4\ \left|R_{e_ie_j}\uu\right|^2, & i,j=1,\ldots,m\\
\tilde K^A(E_i,E_{m+1})=0, & i=1,\ldots,m\\
\tilde K^A(E_i,E_{m+k})=\frac 14\left|R_{\uu e_k}e_i\right|^2, & i=1,\ldots,m,\ k=2,\ldots,m\\
\tilde K^A(E_{m+1}E_{m+k})=-\frac{F_2+2tF_3}{a(t)}, & k=2,\ldots,m\\
\tilde K^A(E_{m+k}E_{m+l})=-\frac{F_2}{a(t)} & k,l=2,\ldots,m.
\end{array}
\right.
\end{equation}

Can we have constant sectional curvature $\tilde c$ on $T(M)$?

If this happens, then it must be $0$, so $T(M)$ is flat.
First, one gets easily that $M$ is locally Euclidean.
Then, we should also have $F_2(t)=0$ and $F_3(t)=0$ for any $t$.
It follows $\MM=\frac \LL{1+2t\LL}$ and $\NN=\frac{\LL'-\LL^2}{1+2t\LL}$. (Hence
$F_1(t)=0$.) These equalities
yield two ordinary differential equations (involving $a$ and $b$), namely:

\quad $(\diamond)\ t(a')^2+2aa'-2ab=0$

\quad $(\diamond\diamond)\ \displaystyle \frac{ab'-2a'b}{a+2tb}=\frac{2a''a-3(a')^2}{2(a+ta')}$.

A simple computation shows that $(\diamond\diamond)$ is consequence of $(\diamond)$. So, we must have
\begin{equation}
b(t)=a'(t)\left(1+\frac{ta'(t)}{2a(t)}\right)\ .
\end{equation}

It is interesting to point our attention to two special cases:

{\bf Case (i)}: $b(t)=\kk a'(t)$, where $\kk$ is a real constant.

\quad \qquad If $a'=0$ then $b=0$ and $a$ is constant, so, $g_A$ is homothetic to Sasaki metric.

\quad \qquad If $a'\neq0$ then $a(t)=a_0t^{2(\kk-1)}$, ($\kk>1$ or $\kk\leq0$, $a_0>0$) and in this case we have to consider
$T_0(M)$.

{\bf Case (ii)}: $b(t)=a(t)$.
We obtain $\frac {a'}a=\frac{-1\pm\sqrt{1+2t}}{t}$ which gives

\qquad \qquad $\displaystyle a(t)=a_0 \frac{e^{2\sqrt{1+2t}}}{(1+\sqrt{1+2t})^2}$, $a_0>0$ \hfill{$(*)$}

or,

\qquad \qquad $\displaystyle a(t)=a_0 \frac{e^{-2\sqrt{1+2t}}}{1+t-\sqrt{1+2t}}$ and in this case we have
to deal with non zero vectors.

\begin{remark}\rm
The manifold $T(M)$ equipped with the Cheeger Gromoll has non constant sectional curvature.
\end{remark}

Putting $a_0=1$ in $(*)$, we can state the following
\begin{proposition}
\label{prop:2.17}
Consider
$g_1$ on $T(M)$ given by
\begin{equation}
\label{metric:g1}
\left\{\begin{array}{l}
g_1(X^H,Y^H)=g(X,Y),\ g_1(X^H,Y^V)=0\\
g_1(X^V,Y^V)=\frac{e^{2\sqrt{1+2t}}}{(1+\sqrt{1+2t})^2}\ \left(g(X,Y)+g(X,\uu)g(Y,\uu)\right)
\end{array}\right.
\end{equation}

The manifold $(T(M),g_1)$ is flat.
\end{proposition}
\proof
For the metric $g_1$ given above we have
$$
\LL=\frac1{1+\sqrt{1+2t}}\ , \MM=\frac1{1+2t+\sqrt{1+2t}}\ ,\ \NN=\frac 1{(1+2t)(1+\sqrt{1+2t})}
$$
$$
F_1=0, \ F_2=0,\ F_3=0.
$$
Thus ${}^1\tilde R=0$.
\rm Here ${}^1\tilde R$ is the Riemannian curvature of the metric $g_1$.

\gata

We can now compare the scalar curvatures of $(M,g)$ and $(T(M),g_A)$.

\begin{proposition}
Let $(M,g)$ be a Riemannian manifold and endow the tangent bundle $T(M)$ with the metric $g_A$. Let
${\rm scal}$ and $\widetilde{\rm scal}^A$ be the scalar curvatures of $g$ and $g_A$ respectively.
The following relation holds:
\begin{equation}
\widetilde{\rm scal}^A={\rm scal}+\frac{2-3a}2\ \sum\limits_{i<j}\left|R_{e_ie_j}\uu\right|^2+
    \frac{1-m}a\ \left(m F_2+4t F_3\right),
\end{equation}
where $\{e_i\}_{i=1,\ldots,m}$ is a local orthonormal frame on $T(M)$.
\end{proposition}
\proof
Using the fact ${\rm scal}=\sum\limits_{i\neq j} K(e_i,e_j)$ and the formula
$\sum\limits_{i,j=1}^m\left|R_{e_i\uu}e_j\right|^2=\sum\limits_{i,j=1}^m\left|R_{e_ie_j}\uu\right|^2$
we get the conclusion.

\gata

\begin{corollary} {\rm (see e.g. \cite{kn:MT88}) }
If $(M,g)$ is a Riemannian manifolds and $T(M)$ is its tangent bundle equipped with
the Sasaki metric $g_S$. Then $(T(M),g_S)$ has constant scalar curvature if and only if
the base manifold $M$ is flat. In this case $T(M)$ is also flat.
\end{corollary}

\begin{corollary}
Let $(M(c),g)$ be a space form and equip $T(M)$ with the metric $g_A$. Then
$$
\widetilde{\rm scal}^A=(m-1)\left[
   mc+t(2-3a)c^2-\frac{mF_2+4tF_3}a
            \right]\ .
$$
\end{corollary}

\begin{corollary} {\rm (see e.g. \cite{kn:Sek91}) }
If $(M,g)$ has constant sectional curvature $c$, then its tangent bundle $T(M)$ endowed
with the Chegeer Gromoll metric $g_{CG}$ is not $($curvature$)$ homogeneous.
\end{corollary}

Could we find functions $a$ and $b$ such that $T(M)$ equipped with the metric $g_A$
has constant scalar curvature?

\vspace{1mm}

First of all consider $a=k$ (a positive real constant). After the computations we
obtain that $b(t)$ should satisfy the following ODE.
\begin{equation}
\label{eq:ODE_ak}
{c^2}\  (2-3\  k)\  {k^3}\  t+b(t)\  \big(k\
\big(m+4\  {c^2}\  (2-3\  k)\  k\  {t^2}\big)+
\end{equation}
$$
+2\ t\  \big(-2+m+2\  {c^2}\  (2-3\  k)\  k\  {t^2}\big)
\  b(t)\big)+2\  k\  t\  {b^{\prime }}(t)=constant.
$$
Let us to give some examples:
\begin{example}\rm
\label{ex:2.22}
If we take $a=\frac 23$ and $b=0$ we obtain that $(T(M),g_A)$ has
constant scalar curvature $\widetilde {scal}^A=m(m-1)c$.
\end{example}
\begin{example}\rm
\label{ex:2.23}
For $a=k=\frac 23$ if the constant in (\ref{eq:ODE_ak}) vanishes, then
(if $b\neq0$) we can integrate the ODE obtaining
$$b=e^{-\frac 32\left[(m-2)t+\frac{m\log t}3\right]}.$$
Then, $(T_0(M),g_A)$ has
constant scalar curvature $\widetilde {scal}^A=m(m-1)c$.
\end{example}
\begin{example}\rm
If we take $a=k\in\big(0,\frac 23\big)$ and $b=\frac{c^2k^2(3k-2)t}{2+m+2c^2(2-3k)kt^2}$ then
$(T(M),g_A)$ has constant scalar curvature $\widetilde {scal}^A=m(m-1)c$.
\end{example}
\begin{example}\rm
If we take $a=1$ and
$b(t)=
\frac{k\  (2+m)+{c^2}\  m\  t}{m\  (2+m)-2\  k\  (2+m)\  t-2\  {c^2}\  m\  {t^2}}
$
then $(T_{t_2}(M),g_A)$ (i.e. the bundle of tangent vectors having length greater than
the positive solution $t_2$ of the equation $m\  (2+m)-2\  k\  (2+m)\  t-2\  {c^2}\  m\  {t^2}=0$)
has constant scalar curvature $\widetilde{scal}^A=(m-1)(mc+k)$, where $k$ is a real constant.
\end{example}

If we consider $b(t)=a(t)$ (as in the case of Cheeger Gromoll metric) then $(T(M),g_A)$
has constant scalar curvature if and only if $a$ satisfies the following ODE:
$$
\begin{array}{l}
-\frac{1}{2\  {{(1+2\  t)}^2}\  {{a(t)}^3}}  \\
\noalign{\vspace{0.5mm}}
\hspace{2.em} \big(-2\  (m+2\  (-2+m)\  t)\  {{a(t)}^2}-4\  t\  {{(c+2\  c\  t)}^2}\  {{a(t)}^3}+  \\
\noalign{\vspace{0.5mm}}
\hspace{4.em} +6\  t\  {{(c+2\  c\  t)}^2}\  {{a(t)}^4}+(-6+m)\  t\  (1+2\  t)\  {{{a^{\prime }}(t)}^2}+
 \\
\noalign{\vspace{0.5mm}}
\hspace{4.em} +2\  a(t)\  ((m+2\  (-1+m)\  t)\  {a^{\prime }}(t)+2\  t\  (1+2\  t)\  {a^\prime{}^\prime}(t))\big)
    ={\rm constant}
\end{array}
$$
which seems to be very complicated to solve.

\section{The tangent spheres bundle}

\subsection{$T_rM$ as hypersurface in $(T(M),g_S,J_S)$}

Let $T_rM=\{v\in T(M)\ :\ g_{\tau(v)}(v,v)=r^2$, with $r\in(0,+\infty)$ and let
$\pi_r:T_rM\longrightarrow M$ be the canonical projection. If we denote by $(x^i,v^i)$
local coordinates on $T(M)$ then, $T_rM$ can be expressed (locally) as
$$g_{\star\star}=r^2,\quad {\rm where} \quad g_{i\star}=g_{ij}v^j$$
($g_{ij}$ are the components of the metric $g$ in the local chart $(U,x^i)$).

Thus, $T_rM$ is a (real) hypersurface in $T(M)$.

We know a generator system for $T_rM$, namely
$\delta_i=\frac\partial{\partial x^i}-\Gamma_{ij}^k(x)v^j\frac\partial{\partial v^k}$
and $Z_i=\frac\partial{\partial v^i}-\frac1{r^2}g_{i\star}v^k\frac\partial{\partial v^k}$,
$i=1,\ldots, m$. (Remark that $\{Z_i\}_{i=1,\ldots,m}$ are not independent
(e.g. $v^iZ_i=0$), but in any point of $T_rM$ they span an $(m-1)$-dimensional
subspace of $TT_rM$.)

Denote by $G_r$ the induced metric from $(T(M), g_S)$; then we can write
\begin{equation}
G_r(\delta_i,\delta_j)=g_{ij},\ G_r(\delta_i,Z_j)=0,\
G_r(Z_i,Z_j)=g_{ij}-\frac1{r^2}\ g_{i\star}g_{j\star}.
\end{equation}

Since the ambient is almost Hermitian, we can construct an almost contact (metric)
structure $(\varphi_r,\xi_r,\eta_r,G_r)$ on $T_rM$ (i.e. $\varphi\in{\mathcal{T}}_1^1(T(M))$,
$\xi_r\in\chi(T(M))$, $\eta_r\in\Lambda^1(T(M))$ verifying $\varphi_r^2=-I+\eta_r\circ\xi_r$,
$\varphi_r\xi_r=0$, $\eta_r\circ\varphi_r=0$, $\eta_r(\xi_r)=1$ and
$G_r(\varphi_rU,\varphi_rV)=G_r(U,V)-\eta_r(U)\eta_r(V)$ (the compatibility condition)).

The unit normal of $T_rM$ is $N_r=\frac1r\ v^i\frac\partial{\partial v^i}$. We put
\begin{equation}
\xi_r=-J_SN_r=\frac 1r\ v^i\delta_i
\end{equation}
(which is unitary and tangent to $T_rM$).

Next, if $U$ is tangent to the hypersurface $T_rM$, define $\varphi_r$ and $\eta_r$ by
$$
\varphi_rU=tan(J_SU)\quad {\rm and} \quad \eta_r(U)N_r=nor(J_SU)
$$
where $tan:TT(M)\longrightarrow TT_rM$ and $nor:TT(M)\longrightarrow N(T_rM)$
are the usual projection operators. (Here $N(T_rM)$ is the normal bundle.)

Hence
\begin{equation}
\left\{
\begin{array}{l}
\varphi_r\delta_i=Z_i, \ \varphi_rZ_i=-\left(\delta_i^h-\frac1{r^2}\ g_{i\star}v^h\right)\delta_h\\
\eta_r(\delta_i)=\frac 1r\ g_{i\star}, \ \eta_r(Z_i)=0.
\end{array}
\right.
\end{equation}

\begin{proposition}
$(T_rM,\varphi_r,\xi_r,\eta_r,G_r)$ is an almost contact metric manifold.
\end{proposition}

We have also
$$
d\eta'_r(\delta_i,\delta_j)=0,\ d\eta'_r(Z_i,Z_j)=0,\
d\eta'_r(\delta_i,Z_j)=-\frac 1{2}\left(g_{ij}-\frac1{r^2}\ g_{i\star}g_{j\star}\right).
$$
(Here we use the formula $d\omega(X,Y)=\frac12\ X\omega(Y)-Y\omega(X)-\omega([X,Y])$, for
$\omega$ a smooth 1-form and for any pair $X,Y$ of vector fields on a manifold.)

In order to have a contact metric structure on $T_r(M)$ (i.e.
$d\eta_r(U,V)=G_r(U,\varphi_r V)$, $\forall$ $U,V\in\chi(T_r(M))$), we have to modify the
almost contact structure in the following way (see e.g. \cite{kn:Bla02}):
$$\varphi^{\rm new}_r=\varphi^{\rm old}_r,\ \xi_r^{\rm new}=2r\xi_r^{\rm old},\
\eta_r^{\rm new}=\frac 1{2r}\ \eta_r^{\rm old} {\ \rm and\ }
G_r^{\rm new}=\frac1{4r^2}\ G_r^{\rm old}.
$$

\subsection{$T_1M$ as hypersurface in $(T(M),g_A,J_A)$}

Let $T_1M=\{\uu\in T(M)\ :\ g_{\tau(\uu)}(\uu,\uu)=1$ and let $\pi:T_1M\longrightarrow M$
be the canonical projection. If we denote by $(x^i, y^i)$ local coordinates on $T(M)$, then
$T_1M$ can be expressed as
$$g_{00}=1,\ {\rm where}\ g_{i0}=g_{ij}y^j.
$$
We know, as above, a generator system for $T_1M$, namely $\delta_i$ and
$Y_i=\frac\partial{\partial y^i}-g_{i0}y^h\frac\partial{\partial y^h}$.
Denote by $G_A$ the induced metric from $(T(M),g_A)$ on $T_1M$; we have
\begin{equation}
G_A(\delta_i,\delta_j)=g_{ij},\ G_A(\delta_i,Y_j)=0,\
G_A(Y_i,Y_j)=a \left(g_{ij} - g_{i0}g_{j0}\right)
\end{equation}
where $a$ is a real positive constant.

Since the ambient manifold is almost Hermitian,
we can construct on $T_1M$ an almost contact metric structure
$(\varphi_A,\xi_A,\eta_A,G_A)$, by using the same method as in previous subsection.
We obtain
\begin{equation}
\left\{
\begin{array}{l}
\varphi_A\delta_i=\frac 1{\sqrt{a}}\ Y_i,\
\varphi_AY_i=-\sqrt{a}\left(\delta_i^h-g_{i0}y^h\right)\delta_h\\
\eta_A(\delta_i)=-\epsilon g_{i0}, \ \eta_A(Y_i)=0, \ \xi_A=-\epsilon y^k\delta_k.
\end{array}
\right.
\end{equation}

We have
\begin{proposition}
$(T_1M,\varphi_A,\xi_A,\eta_A,G_A)$ is an almost contact metric manifold.
\end{proposition}

We have the following expression for the differential $d\eta_A$:
$$d\eta_A(\delta_i,\delta_j)=0, \ d\eta_A(Y_i,Y_j)=0,\ d\eta_A(\delta_i,Y_j)=\frac\epsilon2\ (g_{ij}-g_{i0}g_{j0}).
$$
Similarly to the previous case, in order to have a contact metric structure on $T_1(M)$
we put
$$
\varphi_A^{\rm new}=\varphi_A^{\rm old},\ \xi_A^{\rm new}=-2\epsilon\sqrt{a}\ \xi_A^{\rm old},\
\eta_A^{\rm new}=-\frac \epsilon{2\sqrt{a}}\ \eta_A^{\rm old},\
G_A^{\rm new}=\frac1{4a}\ G_A^{\rm old}.
$$

\subsection{The isometry}

Consider the smooth map
$\tilde F:T(M)\longrightarrow T(M)$ defined by
$\tilde F(p,\uu)=(p,r\uu)$. We will omit in the following
the point $p$.
Remark that if $\uu$ is of unit length, then $r\uu$ has the length $r$, so,
$\tilde F$ restricts to a smooth map
$F:(T_1M,\varphi_A,\xi_A,\eta_A,G_A)\longrightarrow (T_rM,\varphi_r,\xi_r,\eta_r,G_r)$.

We have
\begin{theorem}
The Riemannian manifolds $(T_1M,G_A)$ and $(T_rM,G_r)$ are isometric for $r=\sqrt{a}$.
\end{theorem}
\proof
It is an easy computation to prove that
$dF(\delta_i)=\delta_i$ and $dF(Y_i)=rZ_i$.
Consequently, we have

\qquad $G_r(dF(\delta_i),F(\delta_j))=G_r(\delta_i,\delta_j)=\frac 1{4r^2}\ g_{ij}=\frac a{r^2}\ G_A(\delta_i,\delta_j)$

and

\qquad $G_r(dF(Y_i),dF(Y_j))=G_r(Z_i,Z_j)=\frac 14\ \left(g_{ij}-\frac1{r^2}g_{i\star}g_{j\star}\right)=
\frac 14\ \left(g_{ij}-g_{i0}g_{j0}\right)= G_A(Y_i,Y_j)$.

Hence the conclusion.

\gata

From the contact point of view we can state
\begin{theorem}
$F$ is a $(\varphi_A,\varphi_r)$ map between almost contact manifolds (i.e.
$dF\circ\varphi_A=\varphi_r\circ dF$) if and only if $r=\sqrt{a}$.
\end{theorem}
\proof
One has

\qquad $dF(\varphi_A\delta_i)=-dF(\frac1{\sqrt{a}}\ Y_i)=-\frac r{\sqrt{a}}\ Z_i=
\frac r{\sqrt{a}}\ \varphi_r\delta_i=\frac r{\sqrt{a}}\ \varphi_r dF(\delta_i)$

\qquad $dF(\varphi_A Y_i)=dF(\sqrt{a}\left(\delta_i^h-g_{i0}y^h\right)\delta_h)=
\sqrt{a}\left(\delta_i^h-g_{i0}y^h\right)\delta_h=
\sqrt{a}\left(\delta_i^h-\frac 1{r^2}\ g_{i\star}v^h\right)\delta_h=$

\qquad\qquad\qquad\quad
$=\sqrt{a}\varphi_rZ_i=
\frac{\sqrt{a}}r\ \varphi_r dF(Y_i)$

From here, we get the statement. \gata

\begin{remark} \rm
The characteristic vector field $\xi_A$ is mapped to the characteristic vector field $\xi_r$.
\end{remark}

\subsection{Some properties of $(T_1(M),\varphi_A,\eta_A,\xi_A,G_A)$ as contact manifold}

We have already seen that $(T_1(M),\varphi_A,\eta_A,\xi_A,G_A)$ is a contact manifold.

Denote by $\nabla^A$ the Levi Civita connection on $T_1(M)$ corresponding to the metric $G_A$.
\begin{proposition}
We have
\begin{equation}
\left\{
\begin{array}{l}
\nabla^A_{\delta_i}{\delta_j}=\Gamma^k_{ij}\delta_k-\frac12\ R^k_{0ij}Y_k, \
\nabla^A_{Y_i}\delta_j=\frac a2\ R^k_{j0i}\delta_k\\
\ \\
\nabla^A_{\delta_i}Y_j=\Gamma^k_{ij}Y_k+\frac a2\ R^k_{i0j}\delta_k,\
\nabla^A_{Y_i}Y_j=-g_{j0}Y_i
\end{array}
\right.
\end{equation}
where $R^h_{kij}$ are the local components of the Riemannian curvature on the base manifold $M$
and $"0"$ denotes the contraction with $\uu$, e.g. $R^k_{0ij}=R^k_{lij}y^l$.
\end{proposition}

If we claim that $(T_1(M),\varphi_A,\eta_A,\xi_A,G_A)$ to be a K-contact manifold,
i.e. $\nabla^A_U\xi_A=-\varphi_AU$ for all $U\in\chi(T_1(M))$ (see e.g. \cite{kn:Bla02}),
we can state
\begin{theorem}
The contact metric structure $(\varphi_A,\xi_A,\eta_A,G_A)$ on $T_1(M)$ is K-contact if and
only if the base manifold $(M,g)$ has positive constant sectional curvature $\frac 1a\ $.
In this case $T_1(M)$ becomes a Sasakian manifold.
\end{theorem}
\proof
One can compute
$$
\begin{array}{l}
\nabla^A_{\delta_i}\xi_A=-\sqrt{a}\ R^k_{0i0}Y_k\\
\ \\
\nabla^A_{Y_i}\xi_A=-\sqrt{a}\left[aR^k_{0i0}-2\left(\delta_i^k-g_{i0}y^k\right)\right]\delta_k.
\end{array}
$$
In order to have a K-contact manifold the following relations must occur
\begin{equation}
\label{Kcontact1}
\left\{
\begin{array}{l}
\left(R^k_{0i0}-\frac 1a\delta_i^k\right)Y_k=0\\
\ \\
\left[aR^k_{0i0}-\left(\delta_i^k-g_{i0}y^k\right)\right]\delta_k=0.
\end{array}
\right.
\end{equation}

Since $\{\delta_k\}_{k=1,\ldots,m}$ are linearly independent, it follows
\begin{equation}
\label{Kcontact2}
R^k_{0i0}=\frac 1a\left(\delta_i^k-g_{i0}y^k\right).
\end{equation}
Remark that (\ref{Kcontact2}) implies the first condition in (\ref{Kcontact1}).

We obtain by symmetrization
\begin{equation}
\label{eq:1}
R^k_{lij}+R^k_{jil}=\frac 1a\left(2\delta_i^kg_{jl}-g_{ij}\delta^k_l-g_{il}\delta_j^k\right).
\end{equation}
Using first Bianchi identity one gets
\begin{equation}
\label{eq:2}
2R^k_{lij}+R^k_{ijl}=\frac 1a\left(2\delta_i^kg_{jl}-g_{ij}\delta^k_l-g_{il}\delta_j^k\right).
\end{equation}
Writting now (\ref{eq:1}) after a cyclic permutation and substracting obtained formula from (\ref{eq:2})
one has
\begin{equation}
\label{ew:3}
R^k_{lij}=\frac 1a\left( g_{jl}\delta_i^k-g_{il}\delta_j^k \right)
\end{equation}
which shows that the base manifold $M$ is a real space form $M(\frac 1a)$.

\vspace{2mm}

Conversely, suppose that $M$ has constant sectional curvature $\frac 1a$. Then the curvature can be written
as in (\ref{ew:3}). Let's compute the covariant derivative of $\varphi_A$. We have
$$
\begin{array}{l}
\left(\nabla^A_{\delta_i}\varphi_A\right)\delta_j=\frac 1{2\sqrt{a}}\left(g_{ij}y^h-g_{j0}\delta_i^h\right)\delta_h\\
\ \\
\left(\nabla^A_{\delta_i}\varphi_A\right)Y_j=0\\
\ \\
\left(\nabla^A_{Y_i}\varphi_A\right)\delta_j=-\frac 1{2\sqrt{a}}\ g_{j0}y^i\\
\ \\
\left(\nabla^A_{Y_i}\varphi_A\right)Y_j=\frac{\sqrt{a}}2\left(g_{ij}-g_{i0}g_{j0}\right)y^h\delta_h
\end{array}
$$
which shows that
$$
\left(\nabla^A_U\varphi_A\right)V=G_A(U,V)\xi_A-\eta_A(V)U
$$
for all $U,V\in\chi(T_1(M))$. This relation characterizes Sasakian manifolds among the almost
contact metric manifolds.

This ends the proof.
\gata

\begin{remark}\rm
The tensor field $\varphi_A$ is never parallel.
The manifold $(T_1(M),\varphi_A,\xi_A,\eta_A,G_A)$ is never cosymplectic.
\end{remark}

\vspace{2mm}

{\bf Acknowledgement. } I would like to thank professor Stefano Marchiafava for
inviting me to {\em Istituto di Matematica "G.Castelnuovo", Universit\`a degli
Studi "La Sapienza" Roma}, and for his constant encouragement. I am also grateful
to {\bf CNR} ({\em Consiglio Nazionale delle Ricerche}) Italy, for financial aid.

\

\vspace{2mm}

\scriptsize {

Temporary address:

\sc Universita 'La Sapienza',

Istituto di Matematica 'G. Castelnuovo',

P-le Aldo Moro, n.2,  

00185 - Roma

Italia 

e-mail: \rm munteanu@mat.uniroma1.it,

$\qquad\quad$ munteanu2001@hotmail.com (permanent)

\normalsize

\end{document}